\documentclass[sn-mathphys-num,12pt]{sn-jnl}

\usepackage{amsmath, amsfonts,amsthm,amssymb,amsbsy,upref,color,graphicx,amscd,enumerate}
\usepackage[utf8]{inputenc}

\usepackage[table,xcdraw,svgnames]{xcolor}  

\usepackage{multirow}
\graphicspath{{./metapost/}{./pics/}}

\usepackage{enumitem}

\theoremstyle{plain}
\newtheorem{thm}{Theorem}

\newtheorem{deff}[thm]{Definition}
\theoremstyle{definition}

\newcommand{\ccc}{\operatorname{cyc}}

\newcommand{\bo}[1]{{\bf #1}}

\definecolor{awesome}{rgb}{1.0, 0.13, 0.32}
\definecolor{coolblack}{rgb}{0.0, 0.18, 0.39}
\definecolor{darkcerulean}{rgb}{0.03, 0.27, 0.49}

\makeatletter
\DeclareFontFamily{U}{tipa}{}
\DeclareFontShape{U}{tipa}{m}{n}{<->tipa10}{}
\newcommand{\arc@char}{{\usefont{U}{tipa}{m}{n}\symbol{62}}}%

\newcommand{\arc}[1]{\mathpalette\arc@arc{#1}}

\newcommand{\arc@arc}[2]{%
  \sbox0{$\m@th#1#2$}%
  \vbox{
    \hbox{\resizebox{\wd0}{\height}{\arc@char}}
    \nointerlineskip
    \box0
  }%
}
\makeatother

\title{Optimal Finsler-Hadwiger inequalities}

\author[1,2]{\fnm{Beniamin} \sur{Bogosel}}\email{beniamin.bogosel@polytechnique.edu}

\affil[1]{\orgdiv{Centre de Math\'ematiques Appliqu\'ees}, \orgname{CNRS, \'Ecole Polytechnique}, \orgaddress{\city{Palaiseau}, \postcode{91120}, \country{France}}}

\affil[2]{\orgdiv{Faculty of Exact Sciences}, \orgname{Aurel Vlaicu University of Arad},\newline \orgaddress{\city{2 Elena Dr\u agoi Street}, Arad, \country{Romania}}}

\begin{document}

\abstract{
Various inequalities exist between the area of a triangle, the perimeter squared $(a+b+c)^2$ and the isoperimetric deficit $Q=(a-b)^2+(b-c)^2+(c-a)^2$. The direct and reverse Finsler-Hadwiger inequalities correspond to the best linear inequalities between the three quantities mentioned above. In this paper, the sharpest inequalities between these three quantities are found explicitly. The techniques used involve Blaschke-Santal\'o diagrams and constrained optimization problems.}


\maketitle

\section{Introduction}

Among the most elementary inequalities involving the sides $a,b,c$ of a triangle and its surface area $S$ we have the following:
\begin{equation}\label{eq:weitzenbock}
a^2+b^2+c^2 \geq 4\sqrt{3}S.
\end{equation}
This inequality is attributed to Weitzenb\"ock \cite{Weitzenbock}, but apparently it was proposed more than twenty years earlier by Ionescu as a problem in the Romanian Mathematical Gazette \cite{Ionescu}. Inequality \eqref{eq:weitzenbock} was also proposed as a problem for the third International Mathematical Olympiad in 1961. Nowadays, inequality \eqref{eq:weitzenbock} poses no difficulties to anyone armed with basic knowledge in metric relations in the triangle and basic notions about convexity and Jensen's inequality. A quick proof follows, for example, from the identity:
\begin{equation*}\label{eq:identity} a^2+b^2+c^2 = [(a-b)^2+(b-c)^2+(c-a)^2] + 4S\left(\tan \frac{A}{2}+\tan \frac{B}{2}+\tan \frac{C}{2}\right),
\end{equation*}
where $A,B,C$ denote the measures of the angles of the triangle.

In \cite{Alsina2008} a surprising proof of \eqref{eq:weitzenbock} follows from a geometrical construction involving the Fermat point in a triangle and equilateral triangles constructed on the exterior of each side. The geometric proof actually implies a stronger version of \eqref{eq:weitzenbock} attributed to Finsler and Hadwiger \cite{Finsler1937}:
\begin{equation}\label{eq:HF}
a^2+b^2+c^2 \geq (a-b)^2+(b-c)^2+(c-a)^2 + 4\sqrt{3}S.
\end{equation}
Note that a non negative term was added in the right hand side of \eqref{eq:weitzenbock}. Interestingly, it was observed in \cite{Lukarevski2020} that even though \eqref{eq:HF} is stronger than \eqref{eq:weitzenbock}, applying the Ionescu-Weitzenb\"ock inequality \eqref{eq:weitzenbock} to the triangle made by the midpoints of the arcs determined by the sides on the circumcircle gives precisely the Finsler-Hadwiger inequality. Thus, the Ionescu-Weitzenb\"ock inequality is an example of a \emph{self improving} inequality. 

The following reverse Finsler-Hadwiger inequality also holds 
\begin{equation}\label{eq:reverseHF}
a^2+b^2+c^2 \leq 3[(a-b)^2+(b-c)^2+(c-a)^2] + 4\sqrt{3}S.
\end{equation}
An elementary investigation, looking at triangles with sides $(1,1,\varepsilon), (\frac{1}{2}+\varepsilon,\frac{1}{2}+\varepsilon,1)$ for $\varepsilon\to 0$, shows that the constants before the term $(a-b)^2+(b-c)^2+(c-a)^2$ cannot be improved in \eqref{eq:HF} and \eqref{eq:reverseHF}. Various refinements of the direct and reverse Finsler-Hadwiger (abbreviated FH in the following) inequalities are proposed in \cite{Kouba2017}. 

An improved estimate for the reverse Finsler-Hadwiger inequality is proposed in \cite{LupuEtAl10} for the particular case of acute triangles. The optimal constant for the isoperimetric deficit becomes $(6-\sqrt{6})/2$ in this case:
\[ a^2+b^2+c^2 \leq \frac{6-\sqrt{6}}{2}[(a-b)^2+(b-c)^2+(c-a)^2] + 4\sqrt{3}S.\]

Let us recall the equality cases for inequalities \eqref{eq:HF} and \eqref{eq:reverseHF}:
\begin{itemize}[noitemsep,topsep=0pt]
	\item The direct FH inequality \eqref{eq:HF}: equality holds for equilateral triangles or degenerate triangles with one side of zero length (for example $a=b$, $c=0$).
	\item The reverse FH inequality \eqref{eq:reverseHF}: equality holds for equilateral triangles or degenerate isosceles triangles where one side equals the sum of two others (for example $b=c=a/2$).
\end{itemize}
This observation shows that neither one of the inequalities \eqref{eq:HF}, \eqref{eq:reverseHF} is sharp for general triangles, not verifying the equality cases.

In the following let us denote by 
\begin{equation}\label{eq:Q}
Q:=(a-b)^2+(b-c)^2+(c-a)^2
\end{equation}
the non-negative quantity appearing in the Finsler-Hadwiger inequalities. The notation \eqref{eq:Q} coincides with the one used in the original paper \cite{Finsler1937}. The expression $Q$ defined above can also be called \emph{isoperimetric deficit} since it is a measure of the distance of a triangle with side lengths $a,b,c$ from an equilateral one. 

To make connections with geometrical quantities more obvious, let us reformulate the HF inequalities so that they involve the perimeter. Quick computations show that the Finsler-Hadwiger inequality \eqref{eq:HF} is equivalent to 
\[ 2(ab+bc+ca) \geq a^2+b^2+c^2+4\sqrt{3}S,\]
which, in turn, is equivalent to 
\begin{equation}\label{eq:perim-down}
 (a+b+c)^2 \geq 2Q+12\sqrt{3}S.
 \end{equation}
Furthermore, the reverse Finsler-Hadwiger inequality \eqref{eq:reverseHF} is equivalent to
\[ 6(ab+bc+ca) \leq 5(a^2+b^2+c^2)+4\sqrt{3}S,\]
which is equivalent to 
\begin{equation}\label{eq:perim-up}
 (a+b+c)^2\leq 8Q+12\sqrt{3}S.
\end{equation}
Since the HF inequalities are not sharp except from the particular equality cases enumerated above, the same holds for inequalities \eqref{eq:perim-down}, \eqref{eq:perim-up}. It is, thus meaningful to ask what are the optimal versions of these inequalities.

The goal of the paper in the following is to obtain the sharpest, optimal inequalities involving $a+b+c$, $Q$ (defined in \eqref{eq:Q}) and the area $S$ in the following sense:

\begin{deff}\label{def:optim-ineq}
	Finding the optimal inequalities relating $a+b+c$, $Q$ and $S$  corresponds to finding two functions $F_-, F_+$ which verify
	\[ F_-(a+b+c,Q)\leq S \leq F_+(a+b+c,Q).\]
	Moreover, for any admissible set of values of $(a+b+c,Q)$ there exists a triangle with side lengths $(a,b,c)$ attaining these values and having area $S$. The lower and upper bounds in the above inequality should be attained for at least one triangle.
\end{deff}

Achieving the goal established in Definition \ref{def:optim-ineq} would provide the sharpest inequalities involving the perimeter, the area and the \emph{isoperimetric deficit} $Q$ defined in \eqref{eq:Q}. A simple algebraic reformulation would lead to optimal versions of the HF inequalities \eqref{eq:HF}, \eqref{eq:reverseHF}.

In the following, the functions $F_-, F_+$ from Definition \ref{def:optim-ineq} will be found. The systematic study of complete systems of inequalities linking various geometric quantities related to convex shapes was introduced by Blaschke \cite{Blaschke} and Santal\'o \cite{Santalo}. The key idea is to find a meaningful way to represent the image of quantities of interest on a two dimensional diagram. If the geometry and boundaries of this kind of diagram, called Blaschke-Santal\'o diagram in the following, are known, then explicit and optimal inequalities can be written between the quantities involved.

\section{The Blaschke-Santal\'o diagrams and finding optimal inequalities}

Inequalities \eqref{eq:perim-down}, \eqref{eq:perim-up} are homogeneous with respect to scalings of the triangle sides with the same factor. Let us consider the following two scale invariant ratios involving the terms of interest:
\begin{equation}\label{eq:ratios}  X = \frac{Q}{(a+b+c)^2}, Y = \frac{12\sqrt{3} S}{(a+b+c)^2}.
\end{equation}
Inequality \eqref{eq:perim-down} shows that $Y \in [0,1]$. Precise bounds for $X$ will be given below. Note that the quantities $X,Y$ defined in \eqref{eq:ratios} are well defined as long as the triangle is not reduced to a point, i.e. $a+b+c>0$. It is not difficult to observe that the original HF inequality \eqref{eq:HF} and \eqref{eq:perim-down} are equivalent to
\begin{equation}\label{eq:ublinearXY}
Y \leq 1-2X
\end{equation}
and the reverse HF inequality \eqref{eq:reverseHF}, \eqref{eq:perim-up} is equivalent to
\begin{equation}\label{eq:lblinearXY}
Y \geq 1-8X.
\end{equation}
To include all equality cases in the investigations, flat triangles are allowed. Examples are isosceles triangles where one side has zero length and flat triangles, verifying $b+c=a$. 

Let us consider some numerical simulations, to observe how well the previous bounds characterize the image of the mapping $(a,b,c)\mapsto (X,Y)$ (according to \eqref{eq:ratios}) when $a,b,c$ are general sides of a triangle, given by the Blaschke-Santal\'o diagram
\begin{equation}\label{eq:D}
\mathcal D = \{(X,Y) : \text{ there exist triangle sides }a,b,c\text{ for which }X,Y\text{ are given by \eqref{eq:ratios}}\}.
\end{equation}

\begin{figure}
	\centering 
	\includegraphics[height=0.34\textwidth]{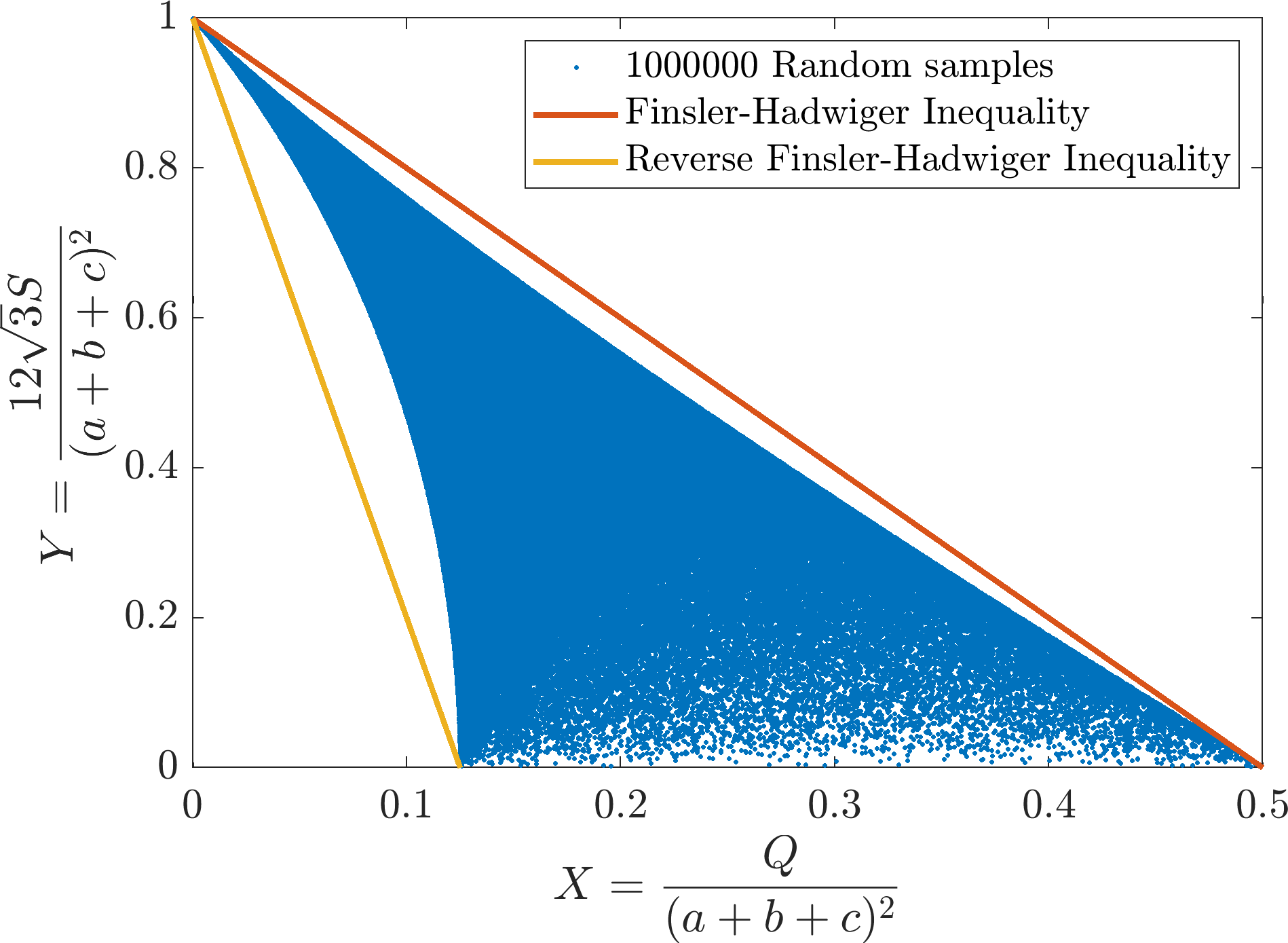}\qquad
	\includegraphics[height=0.34\textwidth]{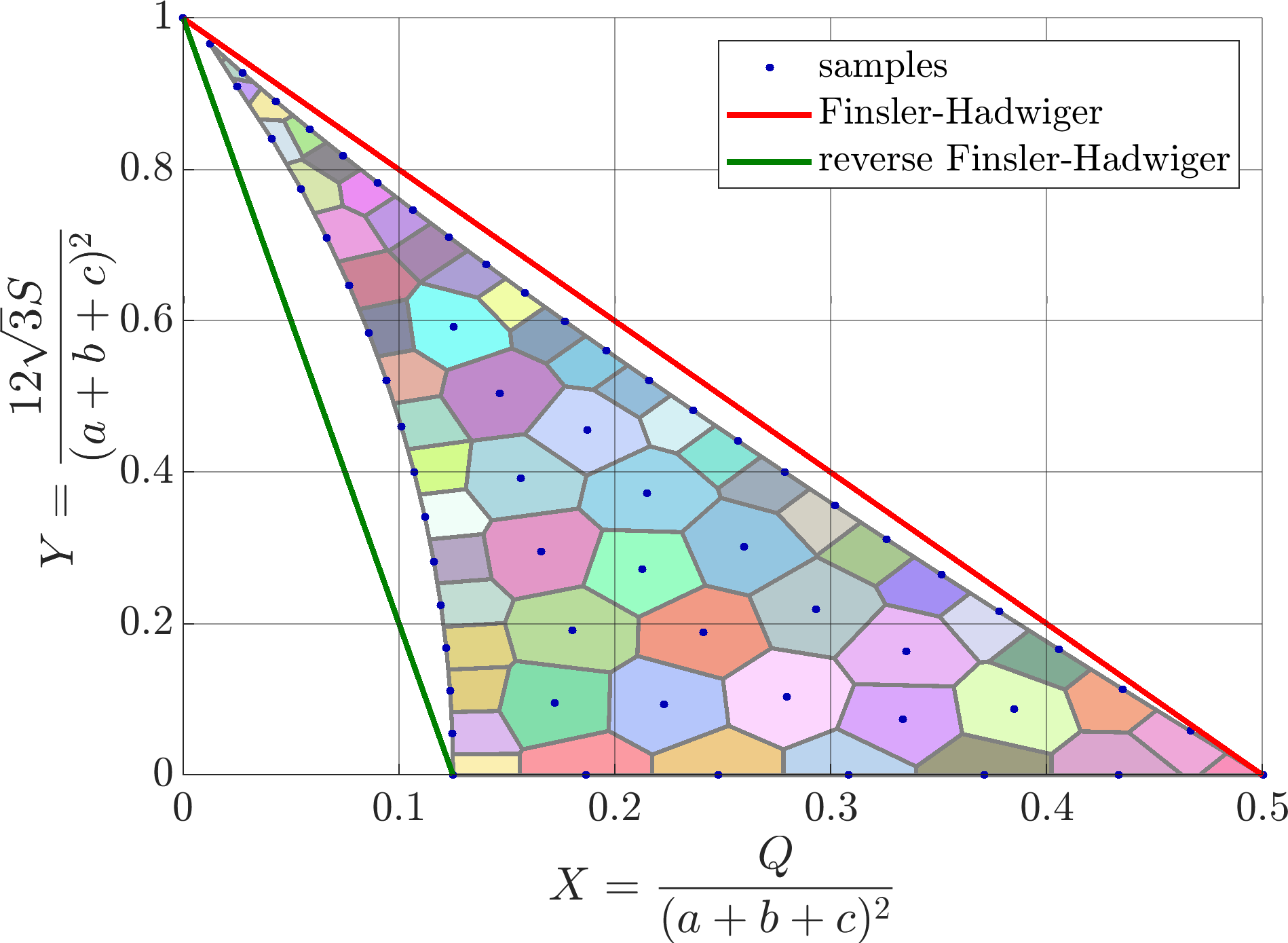}
	\caption{Generating numerically the Blaschke-Santal\'o diagram diagram $\mathcal D$ given by \eqref{eq:D}: random sampling (left) or uniformly distributed samples (right, based on \cite{Blaschke-Santalo}). The direct and reverse Finsler-Hadwiger inequalities correspond to line segments bounding the diagram from above and below.}
	\label{fig:Diagram}
\end{figure}

There exist at least two strategies for approximating the diagram $\mathcal D$ given by \eqref{eq:D}. The first one, easy to implement, consists in generating random triangles, computing $X,Y$ given by \eqref{eq:ratios} and plot the corresponding points. Random triangles can be generated through the classical characterization of the sides of a triangle:
\begin{equation}\label{eq:char} x,y,z\geq 0,  a = y+z, b = z+x, c=x+y.
\end{equation}
Generating $10^6$ random triangles and evaluating $(X,Y)$ from \eqref{eq:ratios} gives the image in Figure \ref{fig:Diagram}.
A second method, described in \cite{Blaschke-Santalo}, proposes to use a limited number of samples whose images are uniformly distributed in the diagram $\mathcal D$. Using only $63$ samples which are well distributed, a more precise visual description of the diagram is found. 

Let us prove some facts observable in Figure \ref{fig:Diagram}. 
\begin{itemize}
	\item $X \leq 1/2$: This is equivalent to $2Q \leq (a+b+c)^2$, or more precisely $a^2+b^2+c^2\leq 2(ab+bc+ca)$. Replacing $a,b,c$ following \eqref{eq:char} gives
	\[ \sum_{\ccc} (x+y)^2 \leq 2\sum_{\ccc} (x+y)(x+z) \Longleftrightarrow \sum_{\ccc}xy \geq 0.\]
	Therefore, $X \leq 1/2$ and equality holds if and only if at least two of $x,y,z$ are zero, implying that one side of the triangle is reduced to a point.
	\item Secondly, the intersection of the diagram with $Y=0$ is given by $X \in [1/8,1/2]$. Indeed, $Y=0$ implies that the triangle has zero area, which means the longest side verifies $a=b+c$. This shows that
	\[ X = \frac{b^2+c^2+(b-c)^2}{4(b+c)^2} = \frac{b^2+c^2-bc}{2b^2+2c^2+4bc}.\]
	A quick computation shows that
	\[ X\leq \frac{1}{2} \Leftrightarrow bc\geq 0 \text{ and } X\geq \frac{1}{8} \Leftrightarrow (b-c)^2\geq 0.\]
	Therefore the desired inequality $X \in [1/8,1/2]$ holds and equality cases correspond to triangles with a side reduced to a point ($X=1/2$) and flat isosceles triangles ($X=1/8$).
	
	Since $X$ is described with a continuous function of $b$ and $c$ it follows that the whole segment $[1/8,1/2]\times \{0\}$ is contained in the image diagram of $(X,Y)$.
	\item The bounds corresponding to the direct and reverse HF inequalities, equivalent to \eqref{eq:ublinearXY} and \eqref{eq:lblinearXY} are represented with line segments in Figure \ref{fig:Diagram}. It was discussed in the introduction that these bounds only touch the diagram at points corresponding to equality cases in the HF inequalities (equilateral triangles, flat isosceles triangles). 
\end{itemize}

In order to find the optimal inequalities, like in Definition \ref{def:optim-ineq}, one needs to characterize the functions giving upper and lower bounds in the diagram shown in Figure \ref{fig:Diagram}. To eliminate the constraint that $a,b,c$ must be the sides of a triangle, the equivalent characterization \eqref{eq:char} is used. This gives, using Heron's formula for the area, 
\begin{equation}\label{eq:ratios-xy} X = \frac{(x-y)^2+(y-z)^2+(z-x)^2}{4(x+y+z)^2}, \ Y = \frac{12\sqrt{3}\sqrt{xyz(x+y+z)}}{4(x+y+z)^2},
\end{equation}
where $x,y,z\geq 0$ are arbitrary.

In the following, we fix the horizontal coordinate $X$ in the diagram $\mathcal D$ and we investigate the possible values that can be taken by the $Y$ coordinate. Since ratios $X,Y$ defined in \eqref{eq:ratios-xy} are scale invariant, without loss of generality, we may assume that $x+y+z=1$. 

Fix $X \in [0,1/2]$, implying that $(x-y)^2+(y-z)^2+(z-x)^2 = 4X$. Therefore 
\[ \left\{\begin{array}{rcl}
2X & = & x^2+y^2+z^2-xy-yz-zx\\
1 & = &x^2+y^2+z^2+2xy+2yz+2zx = (x+y+z)^2
\end{array}\right.\]
Eliminating $x^2+y^2+z^2$, we obtain that 
\[ \frac{1-2X}{3}=xy+yz+zx = xy+z(1-z).\]
Therefore, given $z \in [0,1]$ the numbers $x,y$ verify
\[ x+y = 1-z, \ \ \ xy = \frac{1-2X}{3}-z(1-z),\]
and thus they are roots of the quadratic equation
\[ \lambda^2-(1-z)\lambda+\frac{1-2X}{3}-z(1-z) = 0.\]
The discriminant of this equation is $\Delta(z) = -3z^2+2z-\frac{1}{3}+\frac{8}{3}X$ which is non-negative for 
\[ z \in \left[\frac{1-2\sqrt{2X}}{3},\frac{1+2\sqrt{2X}}{3}\right]=:[z_-,z_+].\]
Observe that since $X \leq 1/2$ the upper bound for $z$ is always at most $1$. Moreover, to obtain $x\geq 0$, $y\geq 0$ one also needs to have $x+y=1-z\geq 0$ and $xy =  \frac{1-2X}{3}-z(1-z)\geq 0$.

Therefore $z\in [0,1]$ and 
\[ z(1-z)\leq \frac{1-2X}{3}.\]
Since $\max_{z \in [0,1]} z(1-z) = 1/4$ we have two situations: 
\begin{itemize}
	\item[(a)] $X \in [0,1/8]$ in which case $(1-2\sqrt{2X})/3\geq 0$ and $(1-2X)/3\geq 1/4$. 
	This implies
	\[ z \in \left[\frac{1-2\sqrt{2X}}{3},\frac{1+2\sqrt{2X}}{3}\right].\]
	\item[(b)] $X \in [1/8,1/2]$ in which case $z(1-z) = (1-2X)/3$ has solutions $\left( 1\pm \sqrt{\frac{8X-1}{3}}\right)/2$. Since $1-2\sqrt{2X}\leq 0$ it follows that
	\[ z \in \left[0,\frac{1-\sqrt{(8X-1)/3}}{2}\right]\cup \left[\frac{1+\sqrt{(8X-1)/3}}{2}, \frac{1+2\sqrt{2X}}{3}\right].\]
\end{itemize}

\begin{figure}
	\centering
	\includegraphics[width=0.45\textwidth]{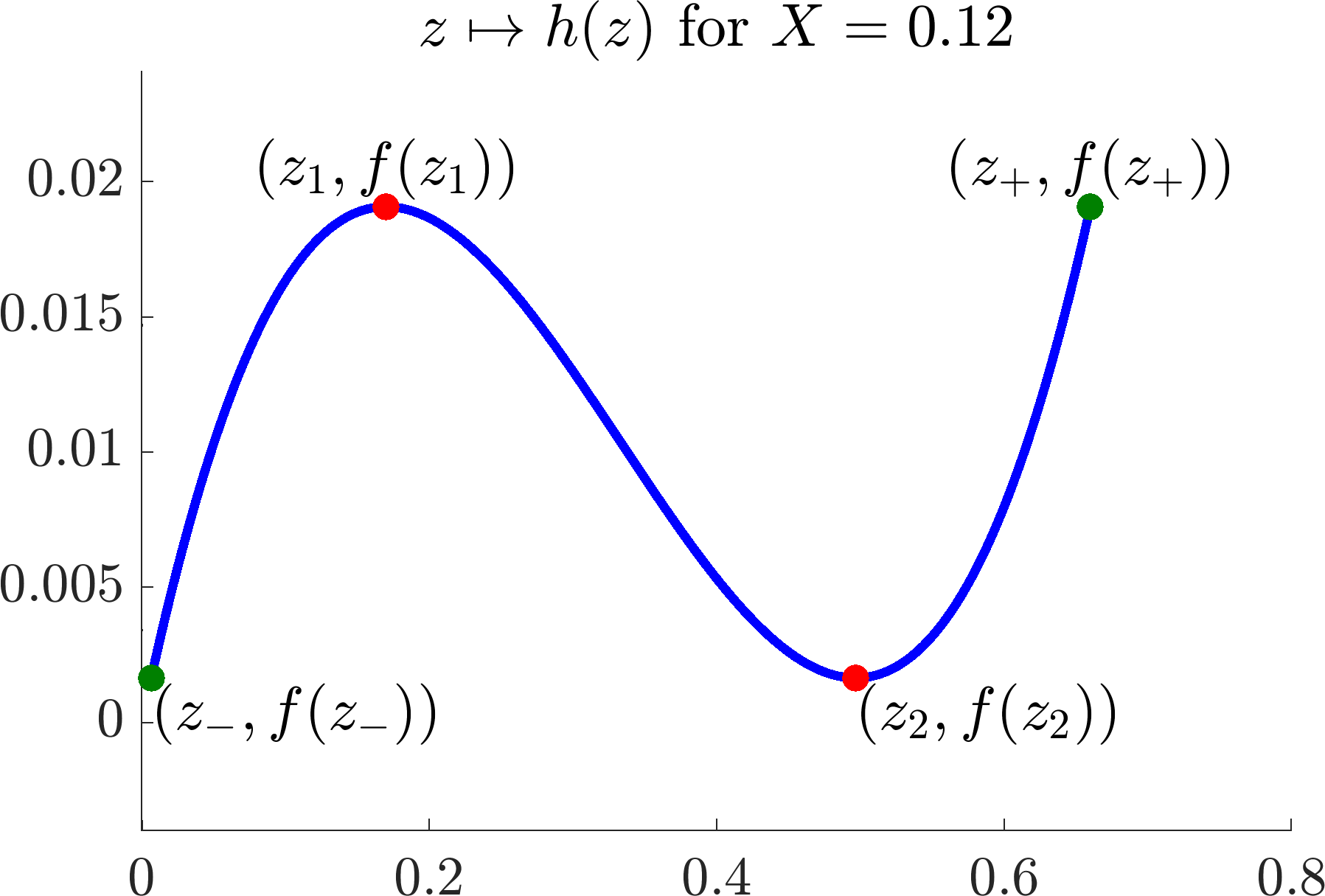}
	\includegraphics[width=0.45\textwidth]{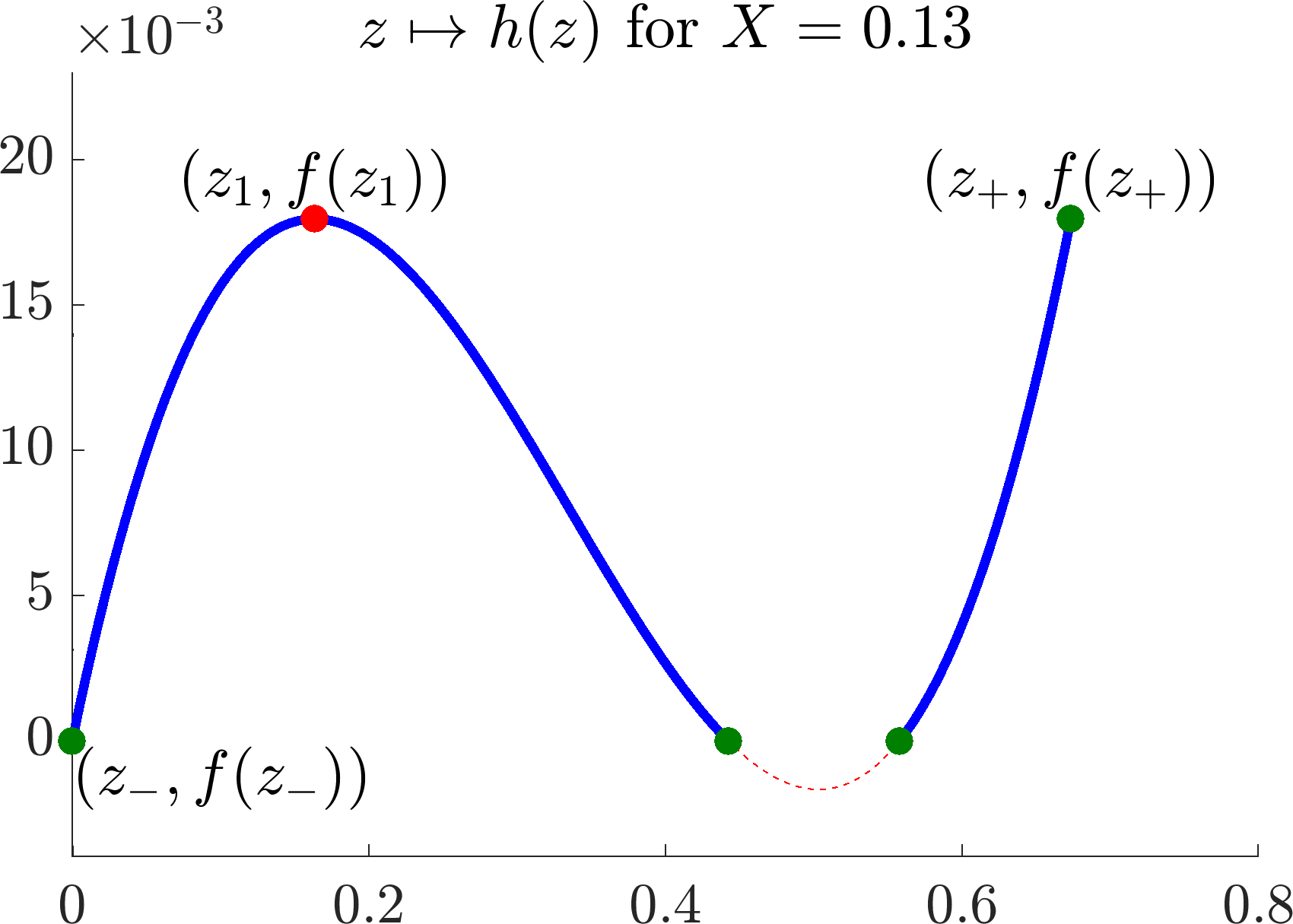}
	\caption{Graph of the function $z \mapsto h(z)$ for admissible values of $z$ in $[z_-,z_+]$ for which $h(z) \geq 0$. (left) Configuration for an example where $X\leq 1/8$ where all values in $[z_-,z_+]$ are admissible. (right) Configuration for $X > 1/8$ where some values of $z$ need to be excluded.}
	\label{fig:graph1D}
\end{figure}

To identify the vertical slices of the diagram one needs to find the image of $Y$ when $X$ is fixed. Under the assumption $x+y+z=1$ we have
\[ Y = 3\sqrt{3}\sqrt{xyz},\]
therefore optimizing $Y$ amounts to optimizing $xyz$, i.e., optimizing
\[ z \mapsto xyz= h(z):=z^3-z^2+\frac{1-2X}{3}z,\]
for $z$ in the range of values which guarantees the existence of $x,y$ with the desired constraints. It remains to investigate the variations of the cubic function $h$, starting with its derivative
\[ h'(z) = 3z^2-2z+\frac{1-2X}{3}.\]
It is immediate to find the zeros of $h'$
\[ z_1 = \frac{1-\sqrt{2X}}{3}, \ \ z_2 = \frac{1+\sqrt{2X}}{3}.\]
Since $h'(z)\leq 0 \Longleftarrow z \in [z_1,z_2]$ we find that
\[ \begin{cases} z \text{ is increasing on } [z_-,z_1]\cup [z_1,z_+] \\
z \text{ is decreasing on } [z_1,z_2].
\end{cases} \]
Graphical representations for the function $z \mapsto h(z)$ are shown in Figure \ref{fig:graph1D} for $X \in \{0.1,0.25\}$, underlining the region where $h$ is non-negative. 

Straightforward computations lead to
\[ \begin{array}{cc}
h\left(\frac{1-2\sqrt{2X}}{3}\right) = h\left(\frac{1+\sqrt{2X}}{3}\right) = \frac{1}{27} \left(-4 \sqrt{2} X^{3/2}-6 X+1\right)\\
h\left(\frac{1+2\sqrt{2X}}{3}\right) = h\left(\frac{1-\sqrt{2X}}{3}\right) = \frac{1}{27} \left(4 \sqrt{2} X^{3/2}-6 X+1\right).
\end{array}\]

Denote 
\[ \phi_-(X)=\sqrt{1-6X-4\sqrt{2}X^{3/2}}, \ \ \phi_-(X)=\sqrt{1-6X+4\sqrt{2}X^{3/2}}.\]

which, in turn, imply the following:
\begin{itemize}
	\item[(a)] If $X \in [0,1/8]$ then, since $Y=3\sqrt{3}\sqrt{h(z)}$ it follows that the image of all possible values of $Y$ is $3\sqrt{3}\sqrt{h([z_-,z_+])}$ which is
	\[ Y \in [\phi_-(X),\phi_+(X)].\]
	\item[(b)] If $X \in [1/8,1/2]$ then the image of all possible values of $Y$ is $3\sqrt{3}\sqrt{h([z_2,z_+])}$ which is
	\[ Y \in [0,\phi_+(X)].\]
\end{itemize}

The computations above lead to the following result: 
\begin{thm}\label{thm:optimal-ineq}
	a) The diagram described in Figure \ref{fig:Diagram} is bounded by: 
	\begin{itemize}[noitemsep,topsep=0pt]
		\item The graph of $X \mapsto \phi_-(X)=\sqrt{1-6X-4\sqrt{2}X^{3/2}}$ on $[0,1/8]$
		\item The graph of $X \mapsto \phi_+(X)=\sqrt{1-6X+4\sqrt{2}X^{3/2}}$ on $[0,1/2]$
		\item The segment $[1/8,1/2]\times \{0\}$.
	\end{itemize}
More precisely:
\[ \phi_-(X)\leq Y \leq \phi_+(X) \text{ for }X \in [0,1/8],\ 0 \leq Y \leq \phi_+(X) \text{ for }X \in [1/8,1/2].\]
The inequalities are sharp or optimal: for a given $X \in [0,1/2]$ all values of $Y$ between the bounds stated above can be attained.

b) Replacing the formulas for $X,Y$ in \eqref{eq:ratios} gives:
\[ 432S^2 \leq  (a+b+c)^4-6Q(a+b+c)^2+4\sqrt{2}Q^{3/2}(a+b+c).\]
If $Q=(a-b)^2+(b-c)^2+(c-a)^2\leq \frac{1}{8}(a+b+c)^2$ then
\[ 432S^2 \geq  (a+b+c)^4-6Q(a+b+c)^2-4\sqrt{2}Q^{3/2}(a+b+c).\]
Both inequalities are sharp or optimal: for a given perimeter $a+b+c$ and isoperimetric deficit $Q=(a-b)^2+(b-c)^2+(c-a)^2$ all values of the area allowed by the above inequalities are attained for some triangle.
\end{thm}

The proof of part a) follows from the previous computations in this section. The proof of part b) comes from replacing the formulas for $X,Y$. It can be observed that inequalities obtained in part b) have the structure required by Definition \ref{def:optim-ineq}.

It is elementary to check that the original Finsler-Hadwiger inequality \eqref{eq:HF} follows from 
\[\sqrt{1-6X+4\sqrt{2}X^{3/2}}\leq 1-2X \text{ for }  X \in [0,1/2]\]
and the reverse one \eqref{eq:reverseHF} follows from
\[\sqrt{1-6X-4\sqrt{2}X^{3/2}}\geq 1-8X \text{ for }  X \in [0,1/8].\]


The method illustrated in this paper provides a general framework for identifying optimal inequalities between three quantities in a triangle:

\bo{Step 1.} Generate two meaningful scale invariant ratios from the three quantities.

\bo{Step 2.} Draw the associated Blaschke-Santal\'o diagram like in Figure \ref{fig:Diagram} using random samples or the more rigorous method shown in \cite{Blaschke-Santalo}. Draw the known inequalities between the quantities of interest to evaluate their optimality.

\bo{Step 3.} Based on the numerical intuition, find theoretically the equations for the curves bounding the diagram. Fixing one of the scale invariant ratios find the optimal bounds for the second one. This can be achieved finding an appropriate parametrization of the variables (like in this article) or using classical methods in constrained optimization theory.

\medskip

\noindent\bo{Declarations:} The authors have no competing interests to declare that are relevant to the content of this article.

\noindent\bo{Data availability:} no data was produced or used while preparing the manuscript.

\bibliographystyle{sn-nature}
\bibliography{./biblio.bib}


\begin{thebibliography}{10}
\ifx \bisbn   \undefined \def \bisbn  #1{ISBN #1}\fi
\ifx \binits  \undefined \def \binits#1{#1}\fi
\ifx \bauthor  \undefined \def \bauthor#1{#1}\fi
\ifx \batitle  \undefined \def \batitle#1{#1}\fi
\ifx \bjtitle  \undefined \def \bjtitle#1{#1}\fi
\ifx \bvolume  \undefined \def \bvolume#1{\textbf{#1}}\fi
\ifx \byear  \undefined \def \byear#1{#1}\fi
\ifx \bissue  \undefined \def \bissue#1{#1}\fi
\ifx \bfpage  \undefined \def \bfpage#1{#1}\fi
\ifx \blpage  \undefined \def \blpage #1{#1}\fi
\ifx \burl  \undefined \def \burl#1{\textsf{#1}}\fi
\ifx \doiurl  \undefined \def \doiurl#1{\url{https://doi.org/#1}}\fi
\ifx \betal  \undefined \def \betal{\textit{et al.}}\fi
\ifx \binstitute  \undefined \def \binstitute#1{#1}\fi
\ifx \binstitutionaled  \undefined \def \binstitutionaled#1{#1}\fi
\ifx \bctitle  \undefined \def \bctitle#1{#1}\fi
\ifx \beditor  \undefined \def \beditor#1{#1}\fi
\ifx \bpublisher  \undefined \def \bpublisher#1{#1}\fi
\ifx \bbtitle  \undefined \def \bbtitle#1{#1}\fi
\ifx \bedition  \undefined \def \bedition#1{#1}\fi
\ifx \bseriesno  \undefined \def \bseriesno#1{#1}\fi
\ifx \blocation  \undefined \def \blocation#1{#1}\fi
\ifx \bsertitle  \undefined \def \bsertitle#1{#1}\fi
\ifx \bsnm \undefined \def \bsnm#1{#1}\fi
\ifx \bsuffix \undefined \def \bsuffix#1{#1}\fi
\ifx \bparticle \undefined \def \bparticle#1{#1}\fi
\ifx \barticle \undefined \def \barticle#1{#1}\fi
\bibcommenthead
\ifx \bconfdate \undefined \def \bconfdate #1{#1}\fi
\ifx \botherref \undefined \def \botherref #1{#1}\fi
\ifx \url \undefined \def \url#1{\textsf{#1}}\fi
\ifx \bchapter \undefined \def \bchapter#1{#1}\fi
\ifx \bbook \undefined \def \bbook#1{#1}\fi
\ifx \bcomment \undefined \def \bcomment#1{#1}\fi
\ifx \oauthor \undefined \def \oauthor#1{#1}\fi
\ifx \citeauthoryear \undefined \def \citeauthoryear#1{#1}\fi
\ifx \endbibitem  \undefined \def \endbibitem {}\fi
\ifx \bconflocation  \undefined \def \bconflocation#1{#1}\fi
\ifx \arxivurl  \undefined \def \arxivurl#1{\textsf{#1}}\fi
\csname PreBibitemsHook\endcsname

\bibitem[\protect\citeauthoryear{Weitzenböck}{1919}]{Weitzenbock}
\begin{barticle}
\bauthor{\bsnm{Weitzenböck}, \binits{R.}}:
\batitle{Uber eine ungleichung in der dreiecksgeometrie}.
\bjtitle{Math. Z.}
\bvolume{5}(\bissue{1-2}),
\bfpage{137}--\blpage{146}
(\byear{1919})
\end{barticle}
\endbibitem

\bibitem[\protect\citeauthoryear{Ionescu}{1897}]{Ionescu}
\begin{barticle}
\bauthor{\bsnm{Ionescu}, \binits{I.}}:
\batitle{Problem 273}.
\bjtitle{Romanian Mathematical Gazette}
\bvolume{13}(\bissue{2}),
\bfpage{52}
(\byear{1897})
\end{barticle}
\endbibitem

\bibitem[\protect\citeauthoryear{Alsina and Nelsen}{2008}]{Alsina2008}
\begin{barticle}
\bauthor{\bsnm{Alsina}, \binits{C.}},
\bauthor{\bsnm{Nelsen}, \binits{R.B.}}:
\batitle{Geometric proofs of the {W}eitzenb\"{o}ck and {H}adwiger-{F}insler
  inequalities}.
\bjtitle{Math. Mag.}
\bvolume{81}(\bissue{3}),
\bfpage{216}--\blpage{219}
(\byear{2008})
\end{barticle}
\endbibitem

\bibitem[\protect\citeauthoryear{Finsler and Hadwiger}{1937}]{Finsler1937}
\begin{barticle}
\bauthor{\bsnm{Finsler}, \binits{P.}},
\bauthor{\bsnm{Hadwiger}, \binits{H.}}:
\batitle{Einige relationen im dreieck}.
\bjtitle{Comment. Math. Helv.}
\bvolume{10}(\bissue{1}),
\bfpage{316}--\blpage{326}
(\byear{1937})
\end{barticle}
\endbibitem

\bibitem[\protect\citeauthoryear{Lukarevski}{2020}]{Lukarevski2020}
\begin{barticle}
\bauthor{\bsnm{Lukarevski}, \binits{M.}}:
\batitle{The circummidarc triangle and the {F}insler-{H}adwiger inequality}.
\bjtitle{Math. Gaz.}
\bvolume{104}(\bissue{560}),
\bfpage{335}--\blpage{338}
(\byear{2020})
\end{barticle}
\endbibitem

\bibitem[\protect\citeauthoryear{Kouba}{2017}]{Kouba2017}
\begin{botherref}
\oauthor{\bsnm{Kouba}, \binits{O.}}:
On certain new refinements of {F}insler-{H}adwiger inequalities.
J. Inequal. Appl.
\textbf{2017}(1)
(2017)
\end{botherref}
\endbibitem

\bibitem[\protect\citeauthoryear{Lupu et~al.}{2010}]{LupuEtAl10}
\begin{botherref}
\oauthor{\bsnm{Lupu}, \binits{C.}},
\oauthor{\bsnm{Mateescu}, \binits{C.}},
\oauthor{\bsnm{Matei}, \binits{V.}},
\oauthor{\bsnm{Opincariu}, \binits{M.}}:
Refinements of the {F}insler-{H}adwiger reverse inequality.
Gazeta {M}atematic\u a ({A}-series)
(28),
49--53
(2010)
\end{botherref}
\endbibitem

\bibitem[\protect\citeauthoryear{Blaschke}{1916}]{Blaschke}
\begin{barticle}
\bauthor{\bsnm{Blaschke}, \binits{W.}}:
\batitle{Eine frage über konvexe köirper}.
\bjtitle{Jahresber. Deutsch. Math. Ver.}
\bvolume{25},
\bfpage{121}--\blpage{125}
(\byear{1916})
\end{barticle}
\endbibitem

\bibitem[\protect\citeauthoryear{Santal\'{o}}{1959/61}]{Santalo}
\begin{barticle}
\bauthor{\bsnm{Santal\'{o}}, \binits{L.A.}}:
\batitle{On complete systems of inequalities between elements of a plane convex
  figure}.
\bjtitle{Math. Notae}
\bvolume{17},
\bfpage{82}--\blpage{104}
(\byear{1959/61})
\end{barticle}
\endbibitem

\bibitem[\protect\citeauthoryear{Bogosel et~al.}{2024}]{Blaschke-Santalo}
\begin{barticle}
\bauthor{\bsnm{Bogosel}, \binits{B.}},
\bauthor{\bsnm{Buttazzo}, \binits{G.}},
\bauthor{\bsnm{Oudet}, \binits{E.}}:
\batitle{On the numerical approximation of {B}laschke–{S}antaló diagrams
  using {C}entroidal {V}oronoi {T}essellations}.
\bjtitle{ESAIM Math. Model. Numer. Anal.}
\bvolume{58}(\bissue{1}),
\bfpage{393}--\blpage{420}
(\byear{2024})
\end{barticle}
\endbibitem

\end{thebibliography}

\bigskip
\small\noindent

\end{document}